\def\Z{\mathbb Z}
  \def\F{\mathbb F}

  \def\a{\alpha}
  
  \def\g{\gamma}
  \def\d{\delta}

  \def\e{\epsilon}

  \def\la{\langle}
  \def\ra{\rangle}
  \def\s{\sigma}

  \def\pf{{\it Proof. }$\;\;$}
  \def\no{\noindent}
  \def\hal{\unskip\nobreak\hfil\penalty50\hskip10pt\hbox{}\nobreak
  \hfill\vrule height 5pt width 6pt depth 1pt\par\vskip 2mm}

    \documentclass[11pt,a4paper]{article}
            \usepackage{amsfonts}
\begin{document}

  \title{Product decompositions in finite simple groups}
  \author{  Martin W. Liebeck \\ Department of Mathematics\\
    Imperial College \\ London SW7 2BZ, UK
    \and
    Nikolay Nikolov \\ Department of Mathematics \\
    Imperial College \\ London SW7 2BZ, UK
    \and
    Aner Shalev \\
    Institute of Mathematics \\
    Hebrew University \\
    Jerusalem 91904, Israel}
  \maketitle

 \begin{abstract}
We propose a general conjecture on decompositions of finite
simple groups as products of conjugates of an arbitrary subset.
We prove this conjecture for bounded subsets of arbitrary finite
simple groups, and for large subsets of groups of Lie type of bounded rank.
Some of our arguments apply recent advances in the theory of growth
in finite simple groups of Lie type, and provide a variety of new
product decompositions of these groups.
 \end{abstract}

  \footnotetext{
  The authors are grateful for the support of an EPSRC grant.
  The third author is also grateful for the support of an ERC
Advanced Grant 247034.}

\newtheorem{theorem}{Theorem}
  \newtheorem{thm}{Theorem}[section]
  \newtheorem{prop}[thm]{Proposition}
  \newtheorem{lem}[thm]{Lemma}
  \newtheorem{cor}[theorem]{Corollary}


In this paper we propose the following conjecture:

\vspace{2mm}
\no {\bf Conjecture } {\it
There exists an absolute constant $c$ such that if
$G$ is a finite simple group and $A$ is any subset of $G$
of size at least two, then $G$ is a
product of $N$ conjugates of $A$ for some $N \le c\log |G|/\log |A|$.}

\vspace{2mm}

Note that we must have $N \ge \log |G|/\log |A|$ by order considerations,
and so the bound above is tight up to a multiplicative constant.

The above conjecture is a stronger version of a recent conjecture we
posed in \cite{LNS}, where $A$ was assumed to be a subgroup of $G$.
Positive evidence for the latter conjecture is provided by \cite{LP} (when
$A$ is a Sylow subgroup) and \cite{LNS1,lub,nik}
(when $A$ is of type $SL_n$), with applications to bounded generation
and expanders.
Further results were proved in \cite{LNS}
in various cases where $A$ is a maximal subgroup of $G$, but the
general case is still open.

In this paper we provide positive evidence for the stronger conjecture
stated above, regarding subsets. One important case where the
conjecture is known to be true (and widely applied) is when the subset $A$
is a conjugacy class, or more generally, a normal subset of $G$; indeed,
this is the main result of \cite{lsann}. Note also that if $G$ is a
group of Lie type of bounded rank, and $A$ is a bounded subset of
$G$, then the conjecture holds, as shown in \cite[2.3]{LNS}.

The following easy reductions will sometimes be useful.
We first claim that, in proving the conjecture for a subset $A$,
we may assume that $1 \in A$.
Indeed, let $a \in A$ and $B = a^{-1}A$.
Then $1 \in B$, and if $G$ is a product of $N$ conjugates of
$B$ then it is also a product of $N$ conjugates of $A$.

Secondly, we claim we may assume there exists $x \ne 1$
such that $1,x,x^{-1} \in A$. Indeed, suppose $1 \in A$
and let $x \in A$ be a non-identity element (whose existence
follows from the assumption $|A| \ge 2$). Then
$1, x, x^2 \in A^2$, hence $x^{-1}, 1, x \in x^{-1}A^2$.
Assuming the conjecture holds for sets containing $x^{-1}, 1, x$
we deduce that $G$ is a product of say
$N \le c \log |G| / \log |A^2| \le c \log |G| / \log |A|$
conjugates of $x^{-1}A^2$, hence it is a product of $N$ conjugates
of $A^2$, so $G$ is a product of $2N \le 2c  \log |G| / \log |A|$
conjugates of $A$.

Our first result here concerns arbitrary subsets of groups of Lie type of
bounded rank, but provides a slightly weaker bound.

\begin{theorem}\label{weak}
Let $G$ be a finite simple group of Lie type of rank $r$, and let
$A$ be any subset of $G$ of size at least $2$. Then there is a constant $c
= c(r)$ depending only on $r$, and a positive integer
$N \le \max (3, (\frac{\log |G|}{\log |A|})^c )$, such that $G$ is a product
of $N$ conjugates of $A$.
\end{theorem}

\pf  The proof is short but relies on strong tools,
most importantly the recent results on growth of Cayley graphs in
\cite{breu,helf,pyber}.
Let $G$ and $A$ be as in the hypothesis.

By \cite[5.3.9]{KL} there exists $\d>0$
depending only on $r$ such that every nontrivial
representation of $G$ has dimension at least $|G|^\d$. Hence \cite[Corollary
1]{NP} shows that if $|A| > |G|^{1-\d/3}$, then $A^3=G$. Consequently
we may assume that $|A| \le |G|^{1-\d/3}$.

Assume (as we may) that $1, x \in A$, where $x \ne 1$.

By \cite[Theorem 2]{HLS}, there are $l = 8(2r+1)$ conjugates
$x^{g_1},\ldots,x^{g_l}$
of $x$ which generate $G$, and hence $G = \la A^{g_1},\ldots,A^{g_l} \ra$.
Define $X = A^{g_1}\cdots A^{g_l}$. Then $X$ contains $A^{g_i}$ for all $i$,
and so $X$ generates $G$. By \cite[Theorem 4]{pyber} or in an equivalent
formulation  \cite[Theorem 2.3]{breu}, for any generating set $Y$ of $G$,
either $Y^3=G$ or $|Y^3|>|Y|^{1+\e}$, where $\e > 0$ depends only on $r$.
(Note that the statements there
only say that $|Y^3|>\gamma |Y|^{1+\e}$ for a positive constant
$\gamma$, but as justified at the beginning of
\cite[Section 6]{helf}, we can
assume $\gamma =1$ by taking a smaller value of $\e$.)
Applying this repeatedly to $X,X^3,X^9,\ldots$, we obtain
\[
|X^{3^n}| \ge {\rm min}(|G|,\,|X|^{(1+\e)^n})
 \ge {\rm min}(|G|,\,|A|^{(1+\e)^n}).
\]
Now choose $n$ minimal such that $(1+\e)^n \ge \frac{\log |G|}{\log |A|}$
and let $k = 3^n$. Then $X^k = G$, and $k \le (\frac{\log |G|}{\log |A|})^b$
where $b = 1+(\log 3/\log (1+\e))$, hence depends only on $r$.
As $X$ is a product of $l$ conjugates of $A$, we see that $G$ is a product
of $kl$ conjugates of $A$. Set $N = kl$. Then $N \le
8(2r+1)(\frac{\log |G|}{\log |A|})^b$. Since $|A| \le |G|^{1-\d/3}$, it follows
that $N \le (\frac{\log |G|}{\log |A|})^c$ for some $c>b$ depending only
on $r$. This completes the proof. \hal
\medskip

Notice that Theorem \ref{weak} implies that if
$A$ is a subset of $G$ of size at least $|G|^\a$ for some fixed
$\a>0$, and $r$ is bounded, then $G$ is a product of boundedly many
conjugates of $A$, so our conjecture holds in this case
(cf. \cite[Theorem 2]{LNS}, which includes an analogous
result for maximal subgroups).

In particular, this leads to a host of new product decompositions,
as follows.

\begin{cor}\label{host}
Let $\bar G$ be a simple adjoint
algebraic group of rank $r$ over the algebraic closure
of $\F_p$, where $p$ is a prime, and let $\s$ be a Frobenius morphism of
$\bar G$ such that $G(q) = (\bar G_\s)'$ is a finite simple group of
Lie type over $\F_q$. Suppose $\bar H$ is a $\s$-stable subgroup of $\bar
G$ of positive dimension such that $H(q) = \bar H_\s \cap G(q)$ is
nontrivial. Then $G(q)$
is equal to a product of $f(r)$ conjugates of $H(q)$, for a suitable function
$f$.
\end{cor}

\pf It is well known that $|H(q)|$ is at least $\g q$ (or $\g q^{1/2}$
for Suzuki and Ree groups), where $\g = \g(r)>0$. As $|H(q)|>1$
by hypothesis, it follows that $|H(q)|>q^\d$ with $\d = \d(r)>0$. Since $|G(q)|
< q^{8r^2}$, the conclusion follows from Theorem \ref{weak}. \hal
\medskip

Various particular cases of this are of special interest.
For example, it follows that a simple group $G(q)$ of rank $r$ (not
a Suzuki group) is a product of $f(r)$ conjugate subgroups
isomorphic to $SL_2(q)$ or $PSL_2(q)$. This was proved in \cite{LNS1,lub}
without the conjugacy part of the conclusion, and was one of the
last steps in showing that all families of finite simple groups can be made
into expanders with respect to  bounded generating sets.

Moreover, various new product decompositions now follow in a similar manner.
For example, $G(q)$ is a product of $f(r)$ conjugates of any nontrivial
torus $T$, or of any centralizer $C_{G(q)}(g)$, and so on.
\medskip

In the second part of this paper we prove our conjecture for
finite simple groups in general, provided the subset $A$ has
bounded size. This follows from the theorem below.

\begin{theorem}\label{strong}
There exists an absolute constant $c$ such that if $G$ is a finite
simple group, and $A$ is any subset of $G$ of size at least $2$,
then $G$ is a product of $N$ conjugates of $A$ for
some $N \le c\log |G|$.
\end{theorem}

\pf Since the case when $G$ is of Lie type and bounded rank follows
from \cite{LNS}, it suffices to prove the theorem for alternating
groups and classical groups of unbounded rank.

We assume (as we may) that $1, x, x^{-1} \in A$ for some $x \ne 1$.

We start with the alternating case $G=A_n$.
It is easy to choose a 3-cycle $y \in A_n$ such that
$[x,y] \ne 1$ has support of size at most 5.
Let $C = x^{A_n}$, the conjugacy class of $x$.
Since $[x,y] = x^{-1}x^y \in C^{-1}C$,
we see that $C^{-1}C$ contains either a 3-cycle,
a 5-cycle or a double transposition. In all cases we deduce
that $(C^{-1}C)^2$ contains all double transpositions in $A_n$.

Since $x, x^{-1} \in A$, some product of 4 conjugates of $A$
contains $\{ 1, t \}$  for a double transposition $t \in A_n$.
Denote by $\tau$ a fixed transposition
of $S_n$, say $(1,2)$. We now use our result that
$S_{n-2}$ (and therefore $A_{n-2}$) is contained in a product
$S_2^{g_1} \cdots S_2^{g_k}$ of $k \leq 320n\log n$ conjugates of
$S_2=\{1, \tau\}$, by Lemma 3.7 with $m=2$ in \cite{LNS}.

By adding the transposition $(n-1,n)$ to the transpositions $\tau^{g_i}$
we obtain a conjugate of $t$ for any copy of $S_2^{g_i}$ and in this way
we see that $A_{n-2}$ is a product of at most $320n \log n$ conjugates
of the set $\{1, t\}$. (We only get even powers of the transposition
$(n-1,n)$ on the last two points since the elements of $A_{n-2}$ always
end up as products of \textbf{even} number of conjugates of $\tau$).

Finally $A_n$ is a product of $3$ conjugates of $A_{n-1}$ (since a product of 2 distinct conjugates of $A_{n-1}$ can move 1 to any point in in $1, \ldots, n$). Therefore $A_{n}$ is a product of $9$ conjugates of $A_{n-2}$. Thus, setting
$c=9 \times 320 \times 4= 11520$, we obtain $A_n$ as a product of
at most $cn \log n$ conjugates of $A$. This concludes the proof for
alternating groups.

Now let $G = PCl_n(q)$,
a projective classical group of (unbounded) dimension $n$ over $\F_q$.
Let $x$ be as above.
By the proof of \cite[2.2]{HLS},
there are elements $y_1,\ldots y_k \in G$ with $k\le 3$ such that the element
$u = [x,y_1,\ldots ,y_k]$ is a non-identity long root element of $G$.
Now $u$ is equal to a product of $2^k \le 8$ conjugates of $x^{\pm 1}$, hence
lies in a product of at most 8 conjugates of $A$.

Replacing $u$ by a conjugate, we may write $u = u_\a(1)$ for a long root
$\a$. Consider the subgroup $H = \langle u_{\pm \a}(t): t \in \F_q \rangle
\cong SL_2(q)$ of $G$. As in \cite[2.3]{LNS}, we see that $H$ is equal to
a product of at most $c_1 \log q$ conjugates of the set $\{1,u\}$, hence it
is contained in a product of at most $c_2 \log q$ conjugates of
$A$.

By \cite{nik}, there is a Levi subgroup $X$ of $G$ of type $SL_m$ such that
$G$ is a product of boundedly many conjugates of $X$; and by
\cite[1.1]{LNS1}, $X$ is a product of at most $c_3n^2$ conjugates of $H$.
Hence $G$ is equal to a product of $c_4n^2$ conjugates of $H$.
We conclude from the above that $G$ is equal to a product of at most
$c_5 n^2 \log q$ conjugates of $A$, completing the proof of the theorem.
\hal

\end{document}